# A REAL SEMINORM WITH SQUARE PROPERTY IS SUBMULTIPLICATIVE

M.El Azhari

Ecole Normale Supérieure, Avenue Oued Akreuch, Takaddoum,

BP 5118, Rabat, Morocco

E-mail: mohammed.elazhari@yahoo.fr

A seminorm with square property on a real associative algebra is submultiplicative.

## 1. INTRODUCTION

It was proved in [3] that every complex seminorm with square property on a commutative algebra is submultiplicative, and it was posed the problem whether this result holds in a noncommutative algebra. This problem was answered in the particular case of Banach algebras [4] and fully resolved in [2], [5] and [7]. The result of [3] holds in the real case. But the results of [2], [4], [5] and [7] don't hold in the real case since we use the Hirschfeld-Zelazko Theorem [6], or its locally bounded version [2], which are not valid in the real case. Using a functional representation theorem [1, Theorem 1], we show that every real seminorm with square property is submultiplicative.

## 2. PRELIMINARIES

Let $A$ be an associative algebra over the field $K = \mathbb{R}$ or $\mathbb{C}$. A seminorm on $A$ is a function $p: A \to [0, \infty)$ satisfying $p(a+b) \leq p(a)+p(b)$ for all $a, b$ in $A$ and $p(ka) = |k|p(a)$ for all $a$ in $A$ and $k$ in $K$. $p$ is a complex seminorm if $K = \mathbb{C}$ and $p$ is a real seminorm if $K = \mathbb{R}$. $p$ is submultiplicative if $p(ab) \leq p(a)p(b)$ for all $a, b$ in $A$. $p$ satisfies the square property if $p(a^2) = p(a)^2$ for all $a$ in $A$. Let $A$ be a real algebra and let $a$ be any element of $A$. The spectrum $sp(a)$ of $a$ in $A$ is defined to be equal to the spectrum of $a$ as an element of the complexification of $A$. If $A$ is unital, then $sp(a) = \{s + it \in \mathbb{C}, (a - se)^2 + t^2 e \notin A^{-1}\}$ for all $a$ in $A$, where $e$ is the unit of $A$ and $A^{-1}$ is the set of all invertible elements in $A$. Let $(A, \|.\|)$ be a real normed algebra, the limit $r(a) = \lim \|a^n\|^{1/n}$ exists for each $a$ in $A$. If $A$ is complete, then $r(a) = \sup\{|v|, v \in sp(a)\}$ for every $a$ in $A$.

## 3. RESULTS

Let (A, $\|.\|$) be a real Banach algebra with unit such that $\|a\| \leq m\, r(a)$ for some positive constant m and all a in A. Let X(A) be the set of all nonzero multiplicative linear functionals from A into the noncommutative algebra H of quaternions. For a in A and x in X(A), put J(a)(x) = x(a). For a in A, J(a): X(A) → H is a map from X(A) into H. X(A) is endowed with the topology generated by J(a), a ∈ A, that is the weakest topology such that all the functions J(a), a ∈ A, are continuous. By [1, Theorem 1], X(A) is a nonempty compact space and the map J: A → C(X(A),H), a → J(a), is an isomorphism (into), where C(X(A),H) is the real algebra of all continuous functions from X(A) into H.

Proposition 3.1. (1) An element a is invertible in A if and only if J(a) is invertible in C(X(A), H), and (2) sp(a) = sp(J(a)) for all a in A.

Proof. (1) The direct implication is obvious. Conversely, there exists g in C(X(A),H) such that J(a) g = g J(a)=1   i.e.   x(a) g(x) = g(x) x(a) = 1  for all x in X(A). Let T be a nonzero irreducible representation of A, by the proof of [1, Theorem 1] there exists S: T(A) → H an isomorphism (into). Since S∘T ∈ X(A) and 0 ≠ S∘T(a) = S(T(a)), it follows that T(a) ≠ 0. If aA ≠ A, there exists a maximal right ideal M containing aA. Let L be the canonical representation of A on A/M which is nonzero and irreducible, also L(a) = 0  since aA is included in M, contradiction. Then aA = A and by the same Aa = A. There exist b, c in A such that ab = ca = e ( e is the unit of A). We have  c = c(ab) = (ca)b = b, so a is invertible in A.

(2) s + it ∈ sp(a)  iff  $(a - se)^2 + t^2 e \notin A^{-1}$

$$\text{iff} \quad J((a-se)^2 + t^2 e) \notin C(X(A),H)^{-1} \quad \text{by (1)}$$

$$\text{iff} \quad (J(a) - s J(e))^2 + t^2 J(e) \notin C(X(A),H)^{-1}$$

$$\text{iff} \quad s + it \in sp(J(a)).$$

Theorem 3.2. Let A be a real associative algebra. Then every seminorm with square property on A is submultiplicative.

Proof. If A is commutative, see [3, Theorem 1]. If A is noncommutative, let p be a seminorm with square property on A. By [5] or [7], there exists m > 0 such that p(ab) ≤ mp(a)p(b) for all a, b in A. Ker(p) is a two sided ideal in A, the norm |.| on the quotient algebra A/Ker(p), defined by |a + Ker(p)| = p(a) is a norm with square property on A/Ker(p). Define $\|a + Ker(p)\| = m |a + Ker(p)|$ for all a in A. Let a, b in A, $\|ab + Ker(p)\| = m|ab + Ker(p)| \leq m^2 |a + Ker(p)| |b + Ker(p)| = \|a + Ker(p)\| \|b + Ker(p)\|$. (A/Ker(p), $\|.\|$) is a real normed algebra. Let a in A, $\|a^2 + Ker(p)\| = m|a^2 + Ker(p)| = m|a + Ker(p)|^2 = m^{-1}(m|a + Ker(p)|)^2 = m^{-1} \|a + Ker(p)\|^2$. The completion B of (A/Ker(p), $\|.\|$) satisfies also the property $\|b^2\| = m^{-1}\|b\|^2$ for all b in B, and consequently $\|b^{2^n}\|^{2^{-n}} = m^{2^{-n}-1}\|b\|$ for all b in B and n in N*, then $r(b) = m^{-1}\|b\|$ i.e. $\|b\| = m\, r(b)$. We consider two cases.

B is unital: By [1, Theorem 1], X(B) is a nonempty compact space and the map J: B → C(X(B),H) is an isomorphism (into). C(X(B),H) is a real Banach algebra with unit under the supnorm $\|.\|_s$. By Proposition 3.1, $r(b) = r(J(b))$ for all b in B. Let b in B, $\|b\| = m\, r(b) = m\, r(J(b)) = m\, \|J(b)\|_s$ since the supnorm satisfies the square property. Then $\|J(b)\|_s = m^{-1} \|b\| = |b|$ for all b in A/Ker(p), so $|.|$ is submultiplicative on A/Ker(p) i.e. p is submultiplicative.

B is not unital: Let $B_1$ be the algebra obtained from B by adjoining the unit. By the same proof of [6, Lemma 2], there exists a norm N on $B_1$ such that

(i)       $(B_1, N)$ is a real Banach algebra with unit
(ii)      $N(b) \leq m^3 \, r_{B_1}(b)$ for all b in $B_1$
(iii)     N and $\|.\|$ are equivalent on B.

By [1, Theorem 1], $X(B_1)$ is a nonempty compact space and the map $J: B_1 \to C(X(B_1), H)$ is an isomorphism (into). Let b in B, $\|b\| = m\, r_B(b) = m\, r_{B_1}(b)$ by (iii)

$$= m\, r(J(b)) \quad \text{by Proposition 3.1}$$

$$= m\, \|J(b)\|_s \quad \text{by the square property of the supnorm.}$$

Then $\|J(b)\|_s = m^{-1} \|b\| = |b|$ for all b in A/Ker(p), so $|.|$ is submultiplicative on A/Ker(p) i.e. p is submultiplicative.